\documentclass[12pt,a4paper]{amsart}

\usepackage{bbm}
\usepackage{lscape}
\usepackage{tabularx,mdwtab}

\usepackage{hyperref}

\usepackage{DLdef1}
\newtheorem*{Assumptions}{\secno Assumptions}

\begin{document}

\title{On reflections of roots in positive characteristic}

\author{Alexey Lebedev}
\address{Equa Simulation AB, Stockholm, Sweden} 
\email{alexey.lebedev@equa.se}

\begin{abstract}
For a Lie superalgebra with Cartan matrix over a field of positive characteristic, some information about its root system in terms of the system of simple roots corresponding to the Chevalley generators is described, under certain given assumptions. This information may be useful for finding other systems of simple roots of such a Lie superalgebra and other Cartan matrices which produce Lie superalgebras isomorphic to a given one.
\end{abstract}

%\date{}

\makeatletter
\@namedef{subjclassname@2020}{\textup{2020} Mathematics Subject Classification}
\makeatother


\maketitle

%\markboth{Reflections}{\itshape Alexey Lebedev}
%{{\itshape Reflections}}

\thispagestyle{empty}

\section{Introduction}

\ssec{Notation and conventions}\label{notation} For a definition of
Lie superalgebras, see  %are defined as in
\cite{KLLS}, especially in characteristics $2$ and $3$.\footnote{For other characteristics, this definition is equivalent to the one usually used in characteristic $0$.} For a definition of Cartan matrices of Lie (super)algebras, the Chevalley generators, roots, and system of simple roots, see %are defined as in
\cite{CCLL}. In what follows, $\fg = \fg(A,I)$ is a Lie superalgebra with Cartan matrix $A$ of rank $n$ and the collection of parities of Chevalley generators $I$, over a~field $\Kee$ of characteristic $p$. Let the elements $e_1^+,\dots , e_n^+, e_1^-,\dots , e_n^-, h_1, \dots, h_n\in\fg$ be the Chevalley generators of $\fg$, let $\alpha_1, \dots, \alpha_n$ be the roots corresponding to $e_1^+,\dots , e_n^+$. For two indices $j$ and $k$ ranging $1$ to $n$ such that $j\neq k$, define $B_{kj}\in \Zee_{\geq 0} \cup \{+\infty\}$ to be
\begin{equation}\label{Bkj}
B_{kj}:=\sup \{m\in \Zee_{\geq 0}\mid \alpha_j + m\alpha_k \text{~is a root}\}.
\end{equation}

I assume that the concept of a Lie superalgebra with Cartan matrix in characteristics $2$ or $3$ is well-defined. Since the definition for Lie superalgebras in these characteristics is not equivalent to the one commonly used in characteristic $0$, this assumption needs to be proven (starting with proving existence of free Lie superalgebras); as far as I know, this has not been done and remains an open problem. %\AL{The part after the semicolon applies to "it needs to be proven", not just to the part in the parentheses}.

\ssec{The aim of the paper} If $\fg$ is a finite-dimensional Lie algebra with Cartan matrix over $\Cee$, and $\alpha_1, \dots, \alpha_n$ is its system of simple roots, reflection in $\alpha_k$ sends this system to another system of simple roots which can be found by the following rule:
\begin{equation}\label{ReflDescr}
\begin{minipage}[c]{12cm} $\alpha_k$ is mapped to $-\alpha_k$; every $\alpha_j$ where $j\neq k$ is mapped to the root of the form $\alpha_j + m\alpha_k$ with the maximal possible $m=B_{kj}$, see \eqref{Bkj}.\end{minipage}
\end{equation}

%As mentioned above, we denote such a maximal $m$ by $B_{kj}$, so the rule can be worded as: $\alpha_k$ is mapped to $-\alpha_k$; every $\alpha_j$, where $j\neq k$, is mapped to $\alpha_j + B_{kj}\alpha_k$.

The same description applies to the transformation of a system of simple roots under an ``odd reflection" (which is not a reflection in the geometrical sense), first described by Serganova in \cite[Appendix]{LSS}.\footnote{Note that in the case of an odd reflection in the root $\alpha_k$, all $B_{kj}$ are either $0$ or $1$.} Given a Cartan matrix and the corresponding Lie superalgebra $\fg = \fg(A,I)$, odd reflections allow one to find other, non-equivalent, Cartan matrices which produce Lie superalgebras isomorphic to $\fg$.

This note is written in the hope that the same approach \eqref{ReflDescr} can be used to find new systems of simple roots with corresponding Cartan matrices whenever $B_{kj}$ are finite for all indices $j\neq k$. The goal of the note is to find explicit expressions for $B_{kj}$ in terms of the Cartan matrix in the case of positive characteristics. It is done under certain assumptions, stated below. These assumptions are known to be true at least in some cases (e.g., for finite-dimensional Lie algebras with Cartan matrix over $\Cee$); but %since
this note covers cases of $p=2$ and $3$, where, as mentioned above, even the existence of Lie superalgebras with Cartan matrices has not been rigorously %strictly
proven, %they are listed as mere
together with the assumptions.

The main result of this note is formula \eqref{main_B_kj}, which describes explicit values of $B_{kj}$ in terms of the Cartan matrix in the case of positive characteristics. A wrong versions of this formula has been published without proof %a derivation
in some previous works, e.g., \cite{BGL}, \cite{BGLL}. %, and others.
This note aims to correct those mistakes.

\section{Main result}

%Let $\fg = \fg(A,I)$ be a Lie superalgebra with Cartan matrix of rank $n$ over field $\Kee$ of characteristic $p>0$. Let $e_1^+,\dots , e_n^+, e_1^-,\dots , e_n^-, h_1, \dots, h_n$ be its Chevalley\DL{} generators, and let $\alpha_1, \dots, \alpha_n$ be the roots corresponding to $e_1^+,\dots , e_n^+$. For two indices $j$ an $k$ ranging\DL{} $1$ to $n$ such that $j\neq k$, define $B_{kj}\in \Zee_{\geq 0} \cup \{+\infty\}$ to be\DL{}$\sup \{m\in \Zee_{\geq 0}\mid \alpha_j + m\alpha_k \text{~is a root}\}$.

In what follows, notations from Subsection \ref{notation} apply. The characteristic $p$ is assumed to be positive. Let $m$ always denote an integer; i.e., even when the upper bound of $m$ is given by $m\leq B_{kj}$, $m$ cannot be $+\infty$ even if $B_{kj}=+\infty$.

%Assume that the following is true:

\begin{Assumptions}\label{assump1-1dim} 1) For all $m$ such that $0\leq m\leq B_{kj}$, the space $\fg_{\alpha_j + m\alpha_k}$ is $1$-dimensional. %\end{Assumption}

%\begin{Assumption}\label{assump2-nonzero}
2) If $m$ is such that $0\leq m<B_{kj}$, then the action of $e_k^+$ on $\fg_{\alpha_j + m\alpha_k}$ is non-zero;

if $m$ is such that $1\leq m\leq B_{kj}$, then the action of $e_k^-$ on $\fg_{\alpha_j + m\alpha_k}$ is non-zero.\end{Assumptions}

\ssbegin{Statement} If the above assumptions are true, then $B_{kj}$ can be found in the following way. Define the sequence $(d_m)_{m=-1}^\infty$ with values in $\Kee$ as follows:

$d_{-1} = 0$;

$d_m = (-1)^{i_k}(d_{m-1} - A_{kj} - mA_{kk})$ for any $m\geq 0$.

Then,  $B_{kj} = \inf \{m\in \Zee_{\geq 0}\mid d_m = 0\}$.\end{Statement}

\begin{proof} For every $m$ such that $0\leq m\leq B_{kj}$, select a non-zero vector $v_m\in\fg_{\alpha_j + m\alpha_k}$. Then,  for any $m$ such that $0\leq m< B_{kj}$, we have $[e_k^+, v_m] = c_m^+ v_{m+1}$ and $[e_k^-, v_{m+1}] = c_m^- v_{m}$ for some non-zero $c_m^+, c_m^-\in\Kee$.

For any $m$ such that $0\leq m< B_{kj}$, denote $d'_m := c_m^+ c_m^-$; also, denote $d'_{-1} := 0$ and, if $B_{kj}<+\infty$, denote $d'_{B_{kj}} := 0$. Then,  for any $m$ such that $0\leq m\leq B_{kj}$, we have $[e_k^+, [e_k^-, v_m]] = d'_{m-1} v_m$ and $[e_k^-, [e_k^+, v_m]] = d'_m v_m$.

Since $[h_k, v_m] = (A_{kj} + mA_{kk})v_m$ is also equal to
\[
[e_k^+, [e_k^-, v_m]] - (-1)^{i_k} [e_k^-, [e_k^+, v_m]] = (d'_{m-1} - (-1)^{i_k} d'_m) v_m,
\]
we see that the  $d'_m$ satisfy the relations for the  $d_m$, i.e., $d'_m = d_m$ for any $m$ such that $0\leq m\leq B_{kj}$. Since $d'_m\neq 0$ for any $m$ such that $0\leq m< B_{kj}$ and $d'_{B_{kj}} = 0$ if $B_{kj}<+\infty$,  the statement is proved.\end{proof}

Now let us find what $B_{kj}$ is equal to, using this approach. First of all, if $A_{kj}=0$, then $d_0=0$, and therefore $B_{kj}=0$.

In what follows, we will assume that $A_{kj}\neq 0$. For any
$x\in \Zee/p\Zee$, we will denote by $x_\Zee$  the minimal non-negative integer congruent to $x$ modulo $p$.

\ssec{$i_k=\ev$} In this case, $d_m = d_{m-1} - A_{kj} - mA_{kk}$ for any $m\geq 0$, and by induction,
\[
d_m = -(m+1)A_{kj} - \binom{m+1}{2}A_{kk}.
\]

\sssec{$i_k=\ev$, $A_{kk}=0$} In this case,
\[
\text{$d_m = -(m+1)A_{kj}$, so $d_m=0$ if and only if $m\equiv -1 \pmod p$,}
\]
 i.e., $B_{kj} = p-1$.

\sssec{$i_k=\ev$, $A_{kk}\neq 0$, $p\neq 2$} In this case, $d_m = -\frac12 (m+1)(2A_{kj} + mA_{kk})$, so $d_{p-1} = 0$, while for any $m$ such that $0\leq m\leq p-2$, we have
\[
\text{$d_m=0$ if and only if $2A_{kj} + mA_{kk} = 0$.}
\]
This means that
$$
B_{kj} = \begin{cases} \left(-\frac{2A_{kj}}{A_{kk}}\right)_\Zee & \text{if $\frac{A_{kj}}{A_{kk}}\in \Zee/p\Zee$;} \\ p-1 & \text{if $\frac{A_{kj}}{A_{kk}}\not\in \Zee/p\Zee$.}\end{cases}
$$

\sssec{$i_k=\ev$, $A_{kk}\neq 0$, $p=2$} In this case,
\[
\begin{array}{ll}
d_0 &= A_{kj}\neq 0;\\
d_1 &= A_{kk}\neq 0;\\
d_2 &= A_{kj}+A_{kk};\\
d_3 &= 0.\\
\end{array}
\]
This means that
 \[
 B_{kj}= \begin{cases} 2&\text{if $A_{kj}=A_{kk}$},\\
  3 &\text{otherwise}.\end{cases}
 \]

\ssec{$i_k=\od$} In this case, $d_m = A_{kj} + mA_{kk} - d_{m-1}$ for any $m\geq 0$. By induction,
$$
d_m = \begin{cases} A_{kj} + lA_{kk} & \text{if $m=2l$, $l\in \Zee_{\geq 0}$;} \\ lA_{kk} & \text{if $m=2l-1$, $l\in \Zee_{\geq 0}$.}\end{cases}
$$

\sssec{$i_k=\od$, $A_{kk}=0$} In this case, $d_0 = A_{kj}\neq 0$, and $d_1=0$, so $B_{kj}=1$.

\sssec{$i_k=\od$, $A_{kk}\neq 0$} In this case, $d_{2p-1} = 0$, and for any $m$ such that $0\leq m \leq 2p-2$, we have $d_m=0$ if and only if $m=2l$, where $l\in \Zee$ and $A_{kj} + lA_{kk}=0$. Therefore,
$$
B_{kj} = \begin{cases} 2\left(-\frac{A_{kj}}{A_{kk}}\right)_\Zee & \text{if $\frac{A_{kj}}{A_{kk}}\in \Zee/p\Zee$;} \\ 2p-1 & \text{if $\frac{A_{kj}}{A_{kk}}\not\in \Zee/p\Zee$.}\end{cases}
$$

Finally, putting this all together, we get:
\begin{equation}\label{main_B_kj}
B_{kj} = \begin{cases} 0 & \text{if $A_{kj}=0$;} \\
p-1 & \text{if $i_k=\ev$, $A_{kk}=0$, $A_{kj}\neq 0$;}\\
& \text{or if $i_k=\ev$, $p\neq 2$, $A_{kk}\neq 0$, $\frac{A_{kj}}{A_{kk}}\not\in \Zee/p\Zee$;}\\
\left(-\frac{2A_{kj}}{A_{kk}}\right)_\Zee & \text{if $i_k=\ev$, $p\neq 2$, $A_{kk}\neq 0$, $\frac{A_{kj}}{A_{kk}}\in \Zee/p\Zee$;}\\
2 & \text{if $i_k=\ev$, $p=2$, $A_{kj}=A_{kk}\neq 0$;}\\
3 & \text{if $i_k=\ev$, $p=2$, $A_{kk}\neq 0$, $A_{kj}\neq A_{kk}$, $A_{kj}\neq 0$;}\\
1 & \text{if $i_k=\od$, $A_{kk}=0$, $A_{kj}\neq 0$;}\\
2\left(-\frac{A_{kj}}{A_{kk}}\right)_\Zee & \text{if $i_k=\od$, $A_{kk}\neq 0$, $\frac{A_{kj}}{A_{kk}}\in \Zee/p\Zee$;} \\
2p-1 & \text{if $i_k=\od$, $A_{kk}\neq 0$, $\frac{A_{kj}}{A_{kk}}\not\in \Zee/p\Zee$.}\end{cases}
\end{equation}
%\AL{Vobshhe govorya, pervyj sluchaj ($A_{kj}=0$) peresekaetsya so sluchayami, gde v uslovii stoIt $\frac{A_{kj}}{A_{kk}}\in \Zee/p\Zee$, no oni dayut odinakovye otvety v sluchayah peresecheniya. StOit li zamenit' na $\frac{A_{kj}}{A_{kk}}\in (\Zee/p\Zee)^\times$, chtoby ne peresekalis'?}\AL{Sluchai, gde otvety 2 i 3, tozhe mozhno zapisat' v terminah togo, lezhit li $\frac{A_{kj}}{A_{kk}}$ v $\Zee/p\Zee$. Mozhet, stOit tak i sdelat'?}
%\DL{ne uveren: tak dlinnee, no zato vse vnyatno skazano}

Note that all the obtained values of $B_{kj}$ for $p>0$ are finite. It is not so for $p=0$: e.g., if $p=0$, $i_k=\ev$, $A_{kk}\neq 0$, $-\frac{A_{kj}}{2A_{kk}}\not\in \Zee_{\geq 0}$, then $B_{kj}=+\infty$.

%\section{Notes}

%\DL{Numbering of Notes is confusing. Why not make them footnotes?}

%\ssec{1} If the assumptions in this text are true, then $B_{kj}$ is always finite for $p>0$. It is not so for $p=0$.

\iffalse\ssec{2} The following assumptions have not been used in these calculations, but they explain their\DL{whose?} purpose:

3) Define $\sigma_1, \dots, \sigma_n$ as follows:
$$
\sigma_j = \begin{cases} \alpha_j + B_{kj}\alpha_k & \text{if $j\neq k$;}\\ -\alpha_k & \text{if $j=k$.}\end{cases}
$$
Then,  $(\sigma_1, \dots, \sigma_n)$ is a system of simple roots in $\fg$.

4) There is a pair $(A', I')$ such that there is an isomorphism between $\fg(A,I)$ and $\fg(A',I')$ which maps $\fh$ to $\fh$ and $\fg(A,I)_{\pm\sigma_i}$ to $\fg(A',I')_{(0, ..., 0, \pm 1, 0, ..., 0)}$, with $\pm 1$  in $i$-th position, for all $i =1, \dots, n$.\fi

\def\eightit{\it}
\def\bib{\bf}
\bibliographystyle{amsalpha}

\begin{thebibliography}{CMZ97}

\bibitem[BGL]{BGL}
Bouarroudj S., Grozman P., Leites D., Classification of finite dimensional modular Lie superalgebras
with indecomposable Cartan matrix.
Symmetry, Integrability and Geometry: Methods and Applications
(SIGMA), 5 (2009), 060, 63 pages; \url{https://arxiv.org/abs/0710.5149}

\bibitem[BGLL]{BGLL}
Bouarroudj S., Grozman P., Lebedev A., Leites D., Divided power (co)homology. Presentations of simple finite dimensional modular Lie superalgebras with Cartan matrix.
Homology, Homotopy and Applications,
Volume 12(1) (2010); \url{https://arxiv.org/abs/0911.0243}


\bibitem[BLLS]{BLLS}
Bouarroudj S., Leites D., Lozhechnyk O., Shang J., The roots of exceptional modular Lie superalgebras with Cartan matrix.  Arnold Math. J. (2020) 6,  63--118;\\
\url{https://arxiv.org/pdf/1904.09578}

\bibitem[CCLL]{CCLL}
Chapovalov~D., Chapovalov~M., Lebedev~A., Leites~D. The
classification of almost affine (hyperbolic) Lie superalgebras.
v.~17, Special issue in memory of F.~Berezin, (2010) 103--161;
\url{https://arxiv.org/abs/0906.1860}


\bibitem[KLLS]{KLLS}
Krutov A., Lebedev A., Leites D., Shchepochkina I., Non-degenerate invariant symmetric bilinear forms on simple Lie superalgebras in characteristic $2$. Oberwolfach preprint OWP 2020-02 \url{http://publications.mfo.de/handle/mfo/3697}

%\bibitem[KLLS1]{KLLS1}
%Krutov A., Lebedev A., Leites D., Shchepochkina I., Non-degenerate invariant symmetric bilinear forms on simple Lie superalgebras in characteristic $2$. Linear Algebra and Its Appl. 649(2) (2022) \url{https://arxiv.org/abs/2102.11653}

\bibitem[LSS]{LSS}
Leites D., Saveliev M., Serganova V., Embeddings of $\mathfrak{osp} (N|2)$ and completely integrable systems. In: M.~Markov,
V.~Man'ko (eds.) \textit{Proc. International Conf. Group-theoretical
Methods in Physics}, %Yurmala, May, 1985. Nauka, Moscow, 1986,
%377--394 MR 89h:17042 (English translation:
VNU Sci Press, (1987), 255--297.

\end{thebibliography}

\end{document}